\documentclass[12]{article}
\usepackage{a4wide}
\usepackage{epsf}
\usepackage{color}
\usepackage{graphicx}
\usepackage{amsmath}
\usepackage{amsthm}
\usepackage{amssymb}
\usepackage{amsfonts}
\usepackage{epsfig}
\usepackage{graphics}
\newtheorem{theorem}{Theorem}[section]

\newtheorem{Remark}{Remark}[section]
\theoremstyle{plain}
\newtheorem{lemma}{Lemma}[section]
\newtheorem{proposition}{Proposition}[section]

\newcommand{\R}{{\mathbb R}}
\newcommand{\C}{{\mathbb C}}

\newcommand{\Z}{{\mathbb Z}}

\newcommand{\im}{{\rm{Im}}}

\newcommand{\I}{{\rm{I}}}
\newcommand{\Ad}{{\rm{Ad}}}
\newcommand{\ad}{{\rm{ad}}}
\usepackage{fancyhdr}

\begin{document}
\begin{titlepage}
\thispagestyle{fancy}
\lhead{Bla}
\pagenumbering{roman}
\author{Maria-Cristina Ciocci\footnote{Research Fellow, Department of
    Mathematics and Statistics,  University of
    Surrey.  Guildford, Surrey, GU2 7XH, UK. E-mail:
    m.ciocci@surrey.ac.uk}, Johan
    Noldus\footnote{Post Doc, Spinoza institute for theoretical physics, Minnaertgebouw,
    Leuvenlaan 4. 3508TD Utrecht, NL.  \emph{Preprint numbers:} SPIN 04/25, ITP-UU 04/42  }  }

\title{A general symmetry preserving reduction scheme and normal form
  for dynamical systems with a compact symmetry group }
\maketitle

\begin{abstract}
We present a generalized Lyapunov Schmidt (\textsc{ls}) reduction
scheme for diffeomorphisms living on a
finite dimensional real vector space $V$ which transform under real
one dimensional characters $\chi$ of an arbitrary compact group with
linear action on $V$.  Moreover we prove a
normal form theorem, such that the normal form still has the desirable
transformation properties with respect
to $\chi$.
\end{abstract}

\end{titlepage}
\pagenumbering{arabic}
\section{Introduction}
We develop a framework to study bifurcation of $q$-periodic points
from a fixed point in families of $\chi$-{\emph{equivariant
maps}}, for a given integer $q\geq 1$. The main question asks for
the solutions of an equation of the form $\Phi_{\lambda}^q(x)=x$
near the fixed point, $x_0$, for $\lambda$ close to a given
critical value $\lambda_0$. Mainly aiming at generalizing the work
of A.\ Vanderbauwhede \cite{Vdb 2,CMC-Vdb 2} and M.C.\ Ciocci
\cite{CMC1} to a broader class of maps, our approach consists of a
combined use of Lyapunov Schmidt ({\sc{ls}}) reduction and normal
form ({\sc{nf}}) reduction. The presented reduction scheme is
structure-preserving and solves the conjectures proposed in
\cite{CMC1}. Note that \cite{CMC1} itself can be seen as a special
case of the theory presented here. The basic idea of the procedure
is as follows. The periodic points we seek are determined by the
zeros of a map $\psi$ on an {\emph{orbit space}} $Y$. By the
Implicit Function Theorem, we reduce $\psi=0$ to
$\psi|_{\widehat{U}}=0$ where $\widehat{U}$ is a subspace of $Y$, larger than the nullspace of the
linearisation of $\psi$ at the bifurcation point. When the map is
put into normal form (in particular, it commutes with the
semisimple part of the linearisation at the bifurcation point) up
to order $k$, the reduced map, and thus the determining equations
for the bifurcations, can be expressed explicitly in terms of this
normal form without going through the details of the reduction. As
we will explain later, the reduced equations retain the full
$\Z_q$-equivariance associated with the {\sc{ls}} reduction and
also inherit the `structure' of the original map. In this
respect, our method is a powerful tool that guarantees as much
symmetry as possible in the reduced equations. One further feature
of the method is that it doesn't require any condition on the
linearisation of the map at the equilibrium (except
invertibility), and can therefore be applied when this
linearisation has purely imaginary eigenvalues of any
multiplicity, and satisfying any number of resonance conditions.
We emphasize that the symmetry induced by the structure of the
original map and the natural $\Z_q$-equivariance associated with
the {\sc{ls}} reduction play an important role in analyzing the
reduced equations and may imply symmetry properties of the
corresponding solutions. For related works and examples we refer for
example to\cite{Lbksymm,LbMelb,CMCPhD}. As our framework applies to problems where one need
solve an equation of the form $P^q(x)=x$ for some map $P$, it
turns out to be particularly useful when investigating  the
existence and bifurcations of subharmonics from a given periodic
orbit in a parameterized family of autonomous systems. Also, the
{\sc{ls}} reduction can be used, in combination  with  ($\Z_q$-)
equivariant singularity theory, to study the
  geometry of resonance tongues obtained by Hopf bifurcation from a
  fixed point of a map. For such a study in the general case (i.e. no
  structure) we refer to \cite{BGV}.

In the following section we specify our setting.
\section{Preliminaries}
We will consider two types of compact groups $(\mathcal{G},.)$
with continuous linear action on a finite dimensional real
vectorspace $V$.  On the one hand $\mathcal{G}$ is a compact
continuous group with left invariant Haar measure $\mu$, and on
the other hand $\mathcal{G}$ is a finite discrete group. Let $\chi
: \mathcal{G} \rightarrow \mathbb{R}$ be a one dimensional
(continuous) real group character \footnote{$\chi$ is irreducible
if the
  corresponding representation is. Moreover $\chi$ satisfies the
  product law $\chi(gh) = \chi(g) \chi(h)$ for all $g,h \in
  \mathcal{G}$ iff the representation is one-dimensional.}, and let
$\psi_{\lambda} : V\rightarrow V$ be a smooth family of local diffeomorphisms
defined on a neighborhood of zero. The map $\psi_{\lambda}$ is called
$\chi${\emph{-equivariant}} if
\begin{equation}  \label{symmetry}
g \circ \psi_{\lambda } \circ g^{-1} = \psi_{ \lambda }^{\chi ( g )}, \quad \forall g \in \mathcal{G}, \end{equation}
with  $\lambda$ belonging to an open neighborhood $\mathcal{O}$ of zero in $\mathbb{R}^{m}$.
Remark that for any $\mathcal{G}$ and one-dimensional character
$\chi$, $\left| \chi (g )\right| = 1$ for any $g \in  \mathcal{G}$.
We shall call $g\in \mathcal{G}$ with $\chi(g)=1$ a \emph{symmetry} of
$\psi$, and $h\in \mathcal{G}$ with $\chi(h)=-1$ a \emph{reversing
  symmetry}.
In the sequel we'll focus on the case $\mathcal{G}$ discrete for
simplicity. But the results admit a straightforward generalization
to the continuous case. The only observation here is that one
should consider averaging over G. That is, $\int$ instead of
$\sum$, (see for example Remark \ref{proj_Gcont}).

\begin{Remark}{\rm{
It is reasonable to think that when a continuous group
$\mathcal{G}$ can be written as $G_0 \rtimes T$, where $G_0$ is
the unit component of $\mathcal{G}$ and $T$ is a finite
subgroup\footnote{This decomposition of $\mathcal{G}$ as a
  semidirect product $\mathcal{G}=G_0\rtimes T$, with $G_0$ the  unit component
  of $\mathcal{G}$ and $T$ a finite subgroup,
is always possible when $\mathcal{G}/G_0$ is abelian,
\cite{BrokD}.}, one could try to gauge out the $G_0$ part by
choosing appropriate coordinates so that only the discrete part
plays a role yielding possible non-trivial real characters. For
example, $O(2)=SO(2)\rtimes \Z_2$ where $0\in \Z_2$ corresponds to
the matrices of det$=1$ and $1$ corresponds to the matrices with
det$=-1$, the product $\rtimes$ is the semidirect product. Any
real one-dimensional character of $O(2)$ is the identity on
$SO(2)$ and $\pm 1$ on $SO(2)\times \{1\}$. Therefore, the
reduction problem, we want to deal with, might reduce to a
discrete problem; i.e., by choosing adapted coordinates one can
gauge out the $SO(2)$ factor and work only with $\Z_2$. Notice
also that the $\Z_2$ factor also restricts the possible complex
one-dimensional representations of the $SO(2)$ factor that one
could choose when the considered family of diffeomorphisms is
generated by a vector field. }}
 \end{Remark}

\medskip\noindent
 We assume that
\begin{equation}\label{basic fixed point}
\psi_{\lambda}(0)=0, \ \ \forall \lambda\in  \mathcal{O} \quad
\mbox{and}
\quad A_{0}:=D_x\psi_0(0) \mbox{ is invertible}.
\end{equation}
\begin{Remark}{\rm{
It should be noted that we do not require that there exists a reversing symmetry that is an involution. }}
\end{Remark}

The goal of this paper is threefold. On the one hand we prove a reduction theorem to
study $q$-periodic points of $\psi_{\lambda}$ in a neighborhood of the
fixed point $0\in \R^n$ for $\lambda \in \mathcal{O}$, (Theorem \ref{reduction}). On the other
hand, using the terminology introduced in \cite{murdock}, we give an
inner product (nilpotent)
normal form for the diffeomorphisms $\psi_{\lambda}$ satisfying the
hypotheses (\ref{symmetry})-(\ref{basic fixed point}), (Theorem \ref{NNF}).
Finally, when the original map is put in normal form, Proposition \ref{relation red diff with NF}  shows
that the reduced map can be approximated without going through the details of the reduction itself.

\section{Reduction}
Let $\psi_{\lambda}:\R^n\times\R^m\rightarrow \R^n$ satisfy the
assumptions given above. Our interest is with the
$q$-periodic points of $\psi_{\lambda}$ in a neighborhood of $0\in
\R^n$ for $\lambda \in \mathcal{O}$. That is,
we look for solutions of the equation
\begin{equation}
\label{ekkie} x = \psi_{\lambda }^{q} (x )
\end{equation}
for $(x , \lambda)$ close to $(0,0) \in V \times \mathbb{R}^{m}$. Denote by $\mathcal{S}_{\lambda}^{q}$ the solution
set, then obviously $\psi_{\lambda}$ defines a $\mathbb{Z}_{q}$ action on $\mathcal{S}_{\lambda}^{q}$.  In order to find
$\mathcal{S}_{\lambda}^{q}$, we lift the equation (\ref{ekkie}) to the
$nq$-dimensional vectorspace $Y_{q}$ of
$q$-periodic sequences, where $n = \textrm{dim} (V )$.
That is: $$Y_{q} = \{ (x_{k})_{k \in  \mathbb{Z}} \, | \, x_{k + q} = x_{k}, \textrm { for all } k \in \mathbb{Z} \}.$$
Let $\widehat{\psi}_{\lambda}$ denote the lift of $\psi_{\lambda}$ to $Y_{q}$, that is
$$ \widehat{\psi}_{\lambda} (( x_{i} )_{i \in \mathbb{Z}}) = ( \psi_{\lambda} ( x_{i} ))_{i \in \mathbb{Z}}$$
and similarly, for all $g \in \mathcal{G}$, let $\widehat{g}$ be
given by
$$ \widehat{g} ((x_{i})_{i \in \mathbb{Z}}) = (g (  x_{\chi ( g )i} ))_{i \in \mathbb{Z}}.$$
Define also the left shift $\sigma$ on $Y_q$ by
$$ \sigma ( (x_{i})_{i \in \mathbb{Z}} ) = (x_{i+1})_{i \in \mathbb{Z}}.$$
Then, equation (\ref{ekkie}) is equivalent to
\begin{equation} \label{2}
\widehat{\psi}_{\lambda} ( x ) = \sigma ( x),
\end{equation}
where $x \equiv (x_{i})_{i \in \mathbb{Z}} \in Y_{q}$. \\*
\\*
\begin{Remark} {\rm{ The following properties are easy to verify.

\begin{itemize}
\item[1-] the `hat' operator, $\ \widehat{} \ $ ,  is a group representation, i.e.,
  $\widehat{g_{1} g_{2}} (y) = \widehat{g_{1}} \widehat{g_{2}} ( y)$ for all $g_{1}, g_{2} \in \mathcal{G}$
and $y \in Y_{q}$,
\item[2-] $\widehat{\psi}_{\lambda}$ is $\sigma$-equivariant, i.e.,
  $\widehat{\psi}_{\lambda} ( \sigma (y)) = \sigma \circ \widehat{\psi}_{\lambda} ( y)$ for all $y \in Y_{q}$,
\item[3-]  $\widehat{g} \circ \widehat{\psi}_{\lambda} \circ
  \widehat{g^{-1}} = \widehat{\psi}_{\lambda}^{\chi (g)} $ for all $g \in \mathcal{G}$,
\item[4-] $\sigma^{q} = id$ and $\widehat{g} \circ \sigma = \sigma^{\chi (g)} \circ \widehat{g}$.
\end{itemize} }}
\end{Remark}
Let $\widehat{A}_{0} = D \widehat{\psi}_{0}(0)$ denote the lift of
$A_0$ to $Y_q$, and let $A_0=N_0+S_0$ be the Jordan-Chevalley
decomposition of $A_0$ (e.g. \cite{CMC1}). Then,  $$\widehat{A}_{0} =
\widehat{S}_{0} + \widehat{N}_{0}$$ is the Jordan-Chevalley decomposition of $\widehat{A}_0$ where
 \begin{itemize}
\item $\sigma  \widehat{S}_{0}=\widehat{S}_0\sigma$ and $\sigma \widehat{N}_{0}=\widehat{N}_0\sigma$,
\item $\widehat{g} \circ \widehat{S}_{0} = \widehat{S}_{0}^{\chi(g)} \circ \widehat{g}$
\end{itemize}

A straightforward application of the (classical) (\textsc{ls})
reduction \cite{Vdb-book} to equation
$\widehat{\Phi}_{\lambda}(y)=\sigma \cdot y$ (\ref{2})  would
result in a bifurcation equation of the form $E(v,\lambda)=0$,
where $E(\cdot, \lambda)$ is a $\Z_q$-equivariant map from
$\ker(D\widehat{\Phi}_0(0)-\sigma)$ into a complement of
$\im(D\widehat{\Phi}_0(0)-\sigma)$ satisfying $E(0,\lambda)=0$ and
$D_{v}E(0,0)=0$. In the case where $A_0=D_x\Phi_0(0)$ and hence
also its lift $\widehat{A}_0=D\widehat{\Phi}_0(0)$ are
non-semisimple, the details of the reduction strongly depend  on
the nilpotent part of $A_0$. Since we do not want to impose any
restriction on $A_0$ except that it has to be inver\-tible, we
perform a \textsc{ls} reduction with respect to the semisimple
part $\widehat{S}_0$ of $\widehat{A}_0$, cf. \cite{CMC1,Vdb
2,CMC-Vdb}.  Recall that `semisimple' means complex diagonizable.

Using the decomposition
\[Y_{q}= \ker\left( \hat{S}_{0}-\sigma \right)\oplus \im\left( \hat{S}_{0}-\sigma \right),\]
as starting point for the reduction of (\ref{2}), we prove the
following result. Note that this theorem proves the conjecture stated
in \cite{CMCPhD} and generalizes \cite{CMC1,CMC-Vdb 2}. The result
also applies to  the reversible-equivariant systems as considered in \cite{Lb locbif}.
\begin{theorem}[Reduction Result]
\label{reduction}
Let $\psi:\R^n\times \R^m \rightarrow \R^n$ be a local family of
$\chi$-equivariant diffeomorphisms,
satisfying (\ref{symmetry}) and (\ref{basic fixed point}).
Let $q\geq 1$, and let $A_0=S_0+N_0$ be the Jordan Chevalley
decomposition of $A_0$. Define the \emph{reduced phase space} $U$ by \[U:=\ker(S_0^q-\I).\] Then there exist a family
of (reduced) diffeomorphisms $\psi_{r,\lambda}:U\rightarrow U$ and a
map $x^*:U\times \R^m \rightarrow\R^n $ such that for each sufficiently small $\lambda \in
\R^m$ the following properties hold:
\begin{description}
\item[(i)] $\psi_{r,\lambda}(0)=0$, $D_u\psi_{r,\lambda=0}=A_0|_{U}$,
  $x^*(0,\lambda)=0$ and, for all $\tilde{u}\in U$, $D_ux^*(0,0)\cdot \tilde{u}=\tilde{u}$ ;
\item[(ii)] $\psi_{r,\lambda}$ is $\Z_q$-equivariant: $\psi_{r,\lambda}(S_0u)=S_0\psi_{r,\lambda}(u)$;
\item[(iii)] $\psi_{r,\lambda}$ is $\chi$-equivariant: $g\circ\psi_{r,\lambda}\circ g^{-1}=\psi_{r,\lambda}^{\chi(g)}$;
\item[(iv)]for sufficiently small $(x,\lambda)\in \R^n\times \R^m$ the
  point $x$ is $q$-periodic under $\psi_{\lambda}$ if and only if
  $x=x^*(u,\lambda)$, where $u\in U$ is $q$-periodic under $\psi_{r,\lambda}$;
\item[(v)] for sufficiently small $(u,\lambda)\in U\times \R^m$ the
  point $u$ is $q$-periodic under $\psi_{r,\lambda}$ if and only if
\begin{equation}\label{intro determining equation I}
\psi_{r,\lambda}(u)=S_0u.
\end{equation}
\end{description}
Moreover, let $\mathcal{B}:U\times \R^n\rightarrow U$  be defined by
\begin{equation}\label{intro reduced bifurcation function}
\mathcal{B}(u,\lambda):=S_0^{-1}\psi_{r,\lambda}(u)-S_0\psi_{r,\lambda}^{-1}(u),
\end{equation}
then
\begin{description}
\item[(vi)] a point $(u,\lambda)\in U\times \R^m$ is a solution of
  equation (\ref{intro determining equation I}) if and only if it is a solution of
\begin{equation}\label{intro bifurcation equation}
\mathcal{B}(u,\lambda)=0;
\end{equation}
\item[(vii)] the map $\mathcal{B}(\cdot,\lambda)$ is such that
\begin{equation}\label{bif eq is Dq equiv}
\mathcal{B}(S_0u,\lambda)=S_0\mathcal{B}(u,\lambda) \quad \mbox{and}
\quad \mathcal{B}(gu,\lambda)=\chi(g)g\mathcal{B}(u,\lambda).
\end{equation}
\end{description}
\end{theorem}
The rest of this section is devoted to the prove of Theorem \ref{reduction}

\begin{lemma} \label{iso}
Let the subspaces $U\subseteq \mathbb{R}^{n}$ and $\widehat{U}\subseteq Y_{q}$ be defined by:
   \begin{equation}
   U:=\ker\left( S_{0}^{q}-I\right) \subseteq \mathbb{R}^{n}, \quad
   \widehat{U}:=\ker\left( \widehat{S}_{0}-\sigma \right) \subseteq Y_{q}.
   \end{equation}
  Then, the mapping
   \begin{equation}
  \xi : U \rightarrow \widehat{U}, \quad \xi \left( u\right) :=\left(S_{0}^iu\right)_{i \in \mathbb{Z}}
   \end{equation}
   is an isomorphism which satisfies the following properties:
\begin{description}
\item[(i)]  $\xi (S_{0}u)=\widehat{S}_{0}\xi (u)=\sigma \cdot \xi (u),$ $u\in U$;
\item[(ii)]  $\xi (A_{0}u)=\widehat{A}_{0}\xi (u)$, $u\in U$;
\item[(iii)] $\widehat{g} \cdot \xi \left( u\right) =\xi \left( gu\right)$, $u\in U$;
\item[(iv)]  $Y_{q}=\xi \left( U\right) \oplus \im\left(
    \widehat{S}_{0}-\sigma \right) $, and this decomposition is invariant under
$\widehat{A}_{0}$, $\widehat{S}_{0}$, and $\sigma $.
\end{description}
\end{lemma}
Note also that $U$ is invariant under $S_{0}$ and $A_{0}$ and that
$\left(\widehat{A}_{0}-\sigma\right)$ is invertible on
$\im(\widehat{S}_{0}-\sigma )$.  The proof of Lemma \ref{iso} is
grosso modo analogous to the proof of Lemma A.2 in \cite{CMC1}.
We should only verify (iii):
   \begin{eqnarray*}
   \widehat{g} \xi(u) &=&\widehat{g} \left( S_{0}^{i}u \right)_{i \in \Z} \\
   & = & \left( gS_{0}^{\chi(g)i}u \right)_{i \in \Z} \\
   &=&\left( S_{0}^{i}gu \right)_{i \in \Z} \\
   & = & \xi(gu).
   \end{eqnarray*}
Obviously, according to item (iv) in the lemma above,
$\sigma \cdot y =\xi(S_0u)+\sigma \cdot v$, where $v \in
\im(\widehat{S}_{0} - \sigma)$ and $u \in U$. It follows that equation (\ref{2}) splits into a system of two equations
   \begin{equation}
   \left\{
   \begin{array}{ll}
   S_{0}u & =\Psi _{\lambda }\left( u,v\right) \text{\quad (a)} \\
   \sigma \cdot v & =\Sigma _{\lambda }\left( u,v\right) \text{\quad (b)}
   \end{array}
   \right.  \label{Reversible system to solveproblem}
   \end{equation}
where the maps
\begin{equation*}
\Psi _{\lambda }:U\times \im\left( \widehat{S}_{0}-\sigma \right)\longrightarrow U \quad \mbox{and} \quad
\Sigma _{\lambda }:U\times \im\left( \widehat{S}_{0}-\sigma \right)\longrightarrow \im\left( \widehat{S}_{0}-\sigma
\right)
\end{equation*}
are uniquely determined by the relation
\begin{equation}\label{for10}
\widehat{\psi }_{\lambda }\left( \xi (u)+v\right) =\xi \left( \Psi _{\lambda}(u,v)\right)
+\Sigma _{\lambda }(u,v).
\end{equation}
One calculates that,
\begin{eqnarray*}
& & \Psi _{\lambda }(0,0)=0, \\
& & \Sigma _{\lambda }(0,0)=0,\\
& &  D_{u}\Psi _{0}(0,0)=\left. A_{0}\right| _{U}, \\
& & D_{v}\Psi _{0}(0,0)=0, \\
& &  D_{u}\Sigma _{0}(0,0)=0, \\
& & D_{v}\Sigma _{0}(0,0)=
\left. \widehat{A}_{0}\right| _{\im\left( \hat{S}_{0}-\sigma \right) }.
\end{eqnarray*}
So, the Implicit Function Theorem applies, hence there exists a unique
mapping $v^{\ast }:U\times \mathbb{R}^{m}\rightarrow \im\left(
   \widehat{S}_{0}-\sigma \right) $, smooth near the
   origin, with $v^{\ast }\left( 0,0\right) =0$ and equation (\ref{Reversible system to
   solveproblem})(b) holds for all
$\left(u,v,\lambda \right) $ $\in U\times \im\left( \widehat{S}_{0}-\sigma
   \right) \times \mathbb{R}^{m}$ if and only if $v=v^{\ast }\left( u,\lambda \right)$.
Observe that the $\mathbb{Z}_{q}$-equivariance of
   $\widehat{\psi }_{\lambda }$ decomposes as
   \begin{equation}
   \Psi _{\lambda }\left( S_{0}u,\sigma \cdot v\right) =S_{0}\Psi _{\lambda}\left( u,v\right),
   \end{equation}
   and
   \begin{equation}
   \Sigma _{\lambda }\left( S_{0}u,\sigma \cdot v\right) =\sigma \cdot \Sigma
   _{\lambda }\left( u,v\right) .
   \label{reversible invariancy property of Big sigma}
   \end{equation}
Moreover, one has that $v^{\ast }\left( 0,\lambda \right) =0$, for
$\lambda \in \mathbb{R}^{m}$ near $0$, and $D_{u}v^{\ast }\left( 0,0\right)
   =0 $. Then, uniqueness of the solution and (\ref{reversible invariancy property
   of Big sigma}) imply that $v^{\ast }\left( S_{0}u,\lambda \right)
 =\sigma \cdot v^{\ast}\left( u,\lambda \right) $.
 \newline
   Substituting the solution $v_{\lambda }^{\ast}\left( u\right) $ into the equation
   (\ref{Reversible system to solveproblem})(a) gives the \emph{determining
   equation}
   \begin{equation}
   S_{0}u=\psi _{r,\lambda }\left( u\right)
   \label{rbe}
   \end{equation}
   where the {\emph{reduced map}}
   $\psi _{r}:U\times \mathbb{R}^{m}\rightarrow U$ is defined by
   \begin{equation}
   \psi _{r,\lambda }\left( u\right) :=\Psi _{\lambda }\left( u,v_{\lambda
   }^{\ast }\left( u\right) \right) .  \label{Reversible reduced map}
   \end{equation}
The following lemma summarizes some basic properties of the reduced map.

\begin{lemma}\label{prm}
The map $\psi_{r, \lambda}$ defined by (\ref{Reversible reduced map}), is such that
\begin{description}
\item[(i)] $\psi _{r,0 }(0)=0$, $\lambda\in \R^m$ and $D_{u}\psi _{r,\lambda }(0)=\left.A_{0}\right| _{U}$;
\item[(ii)] $\psi_{r, \lambda}$ is $\Z_q$-equivariant:  $\psi _{r,\lambda }(S_{0}u)
=S_{0}\psi _{r,\lambda }(u)$, for all $(u,\lambda )\in U\times \mathbb{R}^{m}$.
\item[(iii)] $\psi_{r, \lambda}$ is $\chi$-equivariant: $ g \circ
  \psi_{r, \lambda} \circ g^{-1} =
\psi_{r, \lambda}^{\chi(g)} $.
\end{description}
\end{lemma}

The proof of Lemma \ref{prm} follows from the definitions and the remark that
{\emph{the equation $$\widehat{\psi}_{\lambda}\left( \xi \left( u_{1}\right) +v\right) =
\xi \left(u_{2}\right) +\sigma v$$ holds for all
$u_{1},u_{2}\in U,$ $v\in \im\left( \widehat{S}_{0}-\sigma \right) $ if and only if
$v=v_{\lambda }^{\ast }\left( u_{1}\right)$ and $u_{2}=\psi_{r,\lambda }\left( u_{1}\right). $}
\\* \\*
From here on, the proof of Theorem \ref{reduction} uses the same
arguments as in Theorem 1 of \cite{CMC1} and it is therefore omitted.

Note that:
$$ \widehat{g} v^{*}_{\lambda}(u) =
\sigma^{\frac{1 - \chi(g)}{2}} v^{*}_{\lambda}\Bigl(g \psi_{r,\lambda}^{\frac{1 - \chi(g)}{2}} (u)\Bigr).$$

\subsection{Reduced map via Normal Form}\label{redmapviaNF}

As we will prove later in Section \ref{sec-normalform},  we may assume
that, up to a near-identity transformation, the map  $ \psi_{\lambda }$  has the form
\begin{equation}\label{ReductionNFformula1}
\psi _{\lambda }\left( x\right) =\psi _{\lambda }^{NF}\left( x\right)
+R_{k+1}\left( x,\lambda \right),
\end{equation}
with
\begin{equation*}
 \psi^{NF}_{\lambda}\left( S_{0}x\right) =S_{0}\psi _{\lambda }^{NF}\left(x\right),
\end{equation*}
and
\begin{equation}
R_{k+1}\left( x,\lambda \right) =O\left( \left| x\right| ^{k+1}\right),
\end{equation}
as $x\rightarrow 0$ uniformly in $\lambda $. Now, assuming that
$\psi_{\lambda}$ has been put in the form
(\ref{ReductionNFformula1}), application of the  {\sc{ls}}
reduction scheme as developed above yields that $\psi_{r,\lambda}$
can be approximated as follows.

\begin{proposition} \label{relation red diff with NF}
Assume {\rm{(\ref{basic fixed point})}} and
{\rm{(\ref{ReductionNFformula1})}} then the maps $x^*:U\times \R^m
\rightarrow \R^n$ and $\psi_{r}:U\times \R^m \rightarrow U$ given by Theorem \ref{reduction} are such that
\begin{equation}
x^* (u,\lambda )=u+O(\left| u\right| ^{k+1}) \quad \mbox{and} \quad
\psi _{r,\lambda }(u)=\psi _{\lambda }^{NF}(u)+O(\left| u\right| ^{k+1})
\end{equation}
as $u\rightarrow 0$, uniformly for $\lambda $ in some neighborhood of the
origin of $\mathbb{R}^{m}$. Moreover, $D_{u}\psi _{r}(0,\lambda )=$ $%
D_{x}\left. \psi _{\lambda }^{NF}(0)\right| _{U}=$ $\left.
A_{\lambda }\right| _{U}$ ; so the eigenvalues of $D_{u}\psi
_{r,\lambda }(0)$ coincide with the eigenvalues of $A_{\lambda }$
which are close to $q$'th roots of unity.
\end{proposition}

\begin{Remark} {\rm{In Sec.\ \ref{sec-normalform} we show that things can be arranged
so that $\psi^{NF}$ satisfies some extra constraints involving the
nilpotent part of $D\psi_{0}(0)$. But, this is of no particular
use in determining the reduced diffeomorphism $\psi_{r,\lambda}$.
However, it might be useful when dealing with bifurcation problems
at multiple resonances, see \cite{krein} for an example.}}
\end{Remark}

\section{Linear Normal Form}
Recall from linear algebra (cf.\ e.g.\ \cite{CMC1}) that an operator
$A\in GL(n,\R)$ admits a unique semisimple-unipotent ({\sc{su}})
decomposition $A=S\exp(\mathcal{N})$, where $S$ is semisimple,
$\mathcal{N}$ is nilpotent and $S\mathcal{N}=\mathcal{N}S$.
Goal of this section is to prove that any $A\in GL_{\mathcal{G}}^{\chi}(n,\R)$ close to a given
$A_0\in GL_{\mathcal{G}}^{\chi}(n,\R)$ can be normalised to the form
$A=S_0e^{\mathcal{N}_0+B}$, where $A_0=S_0e^{\mathcal{N}_0}$ is the
SU-decomposition of $A_0$ and $S_0B=BS_0$,
$B\mathcal{N}^T_0=\mathcal{N}_0^TB$. Such normal form depends on the
choice of a suitable inner product for  $gl(\R,n)$, whose existence is
proved first (see Lemma \ref{inprod1} and Lemma \ref{sp} below).

We start with some technical details.

\subsection{Technicalities}
Let $\chi_{1}$, $\chi_{2}$ be one-dimensional real characters of a
compact finite discrete group $\mathcal{G}$.  Define the associated projection operators
$ P^{\chi_{i}}_{\mathcal{G}}: gl(n, \R) \rightarrow gl(n, \R)$ by
\begin{equation}
 P^{\chi_{i}}_{\mathcal{G}}(A) := \frac{1}{\left| \mathcal{G} \right|}
 \sum_{g \in \mathcal{G}} \chi_{i}(g) gAg^{-1}.
\end{equation}

\begin{Remark}\label{proj_Gcont} {\rm{In the continuos case one should consider
\[ P^{\chi_{i}}_{\mathcal{G}}(A) := \int_ \mathcal{G} \chi_{i}(g) gAg^{-1}.\] }}
\end{Remark}

\medskip
\noindent
If $\chi_{1}$ and $\chi_{2}$ are orthogonal characters, i.e.,
$\sum_{g \in \mathcal{G}} \chi_{1}(g^{-1}) \chi_{2}(g) = 0$
then $$ P^{\chi_1} P^{\chi_2} = 0.$$
Define the subset $gl^{\chi}_{\mathcal{G}}(n , \R) \subset gl(n, \R)$ by
\begin{equation}
gl^{\chi}_{\mathcal{G}}(n, \R) :=
P_{\mathcal{G}}^{\chi} (gl(n,\R)) =\left\{ A \in gl(n, \R) |\,
  gAg^{-1} =  {\chi(g)}A \right\}.
\end{equation}
and the set $GL_{\mathcal{G}}^{\chi}(n,\R)$ by
\begin{equation}
 GL^{\chi}_{\mathcal{G}}(n, \R) :=  \left\{ A \in GL(n, \R) |\,
 gAg^{-1} = A^{\chi(g)} \right\}.
\end{equation}

\begin{Remark} {\rm{Note that $GL_{\mathcal{G}}^{\chi}(n,\R)$ is in general only an
algebraic variety and not a manifold. However, this does not occur
in the case of reversible systems as considered in \cite{CMC1}.
For example, $GL_{\mathcal{G}}^{\chi}(n,\R)$ is, in the simple
case $n=2$, $\mathcal{G}=\{\I, R\}$ with $R\in GL(2,\R)$ such that
$R^2=\I$, and $\chi(R)=-1$, diffeomorphism equivalent to a
hyperboloyde.}}
\end{Remark}

\medskip
For $A \in GL(n, \R)$, define the adjoint action $Ad(A)$ of $A$ on $gl(n, \R)$ by
$$ Ad(A)B = ABA^{-1}. $$
We need to study the action of $Ad(A^{-1}_{0}) - I$ on the space
$gl^{\chi}_{\mathcal{G}}(n, \R)$ of all $\chi$-equivariant operators,
where $A_{0} \in GL_{\mathcal{G}}^{\chi}(n, \R)$.
Therefore, we introduce the group $\mathcal{G}^{\chi}(A_{0})$,
$$\mathcal{G}^{\chi}(A_{0}) = \left\{ x | x \in \mathcal{G} \textrm{
    and } \chi(x) = 1, \textrm{ or } x
= g A_{0} \textrm{ where } g \in \mathcal{G}, \chi(g) = -1 \right\}. $$
A straightforward calculation shows that
\begin{eqnarray}
\label{trans} Ad(A_{0}^{-1}) - I : gl^{\chi}_{\mathcal{G}}(n, \R) \rightarrow gl^{1}_{\mathcal{G}^{\chi}(A_{0})}(n,\R).
\end{eqnarray}
Given a character $\alpha$ of the group $\mathcal{G}$, we can define
an associated character $\widetilde{\alpha}$ of the group $\mathcal{G}^{\chi}(A_{0})$ as follows:
\begin{itemize}
\item $\widetilde{\alpha}(h) = \alpha(h)$, if $\chi(h) = 1$ and $h \in \mathcal{G}$;
\item $\widetilde{\alpha}(gA_{0}) = \alpha(g)$, if $\chi(g)=-1$ and $g \in \mathcal{G}$.
\end{itemize}
With this notation, one shows that the following generalization of (\ref{trans}) holds:
$$ Ad(A_{0}^{-1}) - I: gl^{\alpha}_{\mathcal{G}}(n, \R) \rightarrow gl^{\widetilde{\alpha} \widetilde{\chi}}_{G^{\chi}(A_{0})}(n, \R).$$

\subsection{Inner Product}
\begin{lemma} \label{inprod1}
Let $S_{0} \in GL_{\mathcal{G}}^{\chi}(n,\R)$  be semisimple and
$\chi$ be as before, then there exists a scalar product on $\R^{n}$
(with corresponding involution $*$) such that $S_0$ is normal, i.e.,
$S_{0}S_{0}^{*} = S_{0}^{*}S_{0}$, and $g S_{0}^{*} = (S_{0}^{*})^{\chi(g)} g$, for all $g \in \mathcal{G}$.
\end{lemma}
\begin{proof}
Denote by $\sigma(S_0)$ the spectrum of $S_0$ and by $(\cdot,
\cdot)$ the usual scalar product on $\R^n$. Note that $S_{0} \in
GL_{\mathcal{G}}^{\chi}$ implies that the spectrum of $S_0$
consists of $\pm 1$'s, or pairs $\{\lambda,
\overline{\lambda}\}\in \C$ with $|\lambda|=1$, or quadruples
$\{\lambda,
\lambda^{-1},\overline{\lambda},\overline{\lambda}^{-1} \}$. Since
$S_0$ is semisimple, the spectral decomposition holds, i.e.,
$$\C^n=\oplus_{\lambda \in \sigma(S_0)} V_\lambda$$
where $V_{\lambda}$
denotes the eigenspace corresponding to the eigenvalue $\lambda$.
Define a scalar product $\langle \cdot, \cdot \rangle$ on $\R^n$ as follows:
\[\langle x,x^{\prime}\rangle:= \sum_{\lambda \in \sigma(S_0)\cap \R}(v_\lambda,v^\prime_\lambda)_{\lambda}+
\sum_{\lambda \in \sigma(S_0), \im(\lambda)>0}(w_{\lambda},w^{\prime}_{\lambda})_{\lambda},\]
where
$v_{\lambda},v^{\prime}_{\lambda}\in V_{\lambda}\cap\R^n$,
$w_{\lambda}, w^{\prime}_{\lambda}\in
W_{\lambda}:=\left(V_{\lambda}\oplus V_{\overline{\lambda}}\right)\cap
\R^n$, and moreover $x$ decomposes uniquely as $x=\sum_{\lambda
  \in\sigma(S_0)\cap \R}v_{\lambda}+\sum_{\lambda \in \sigma(S_0),
  \im(\lambda)>0}w_{\lambda}$ and similarly $x^{\prime}$ does.
Our aim is to determine $(\cdot, \cdot)_{\lambda}$ such that $S_0$ is normal with respect to this scalar product.

When $\lambda \in \R$, any choice of scalar product is good since $S_0|_{V_{\lambda}}=\lambda\I_{V_{\lambda}}$.
So, the first part of the proof is completed when our assertion is
shown for $S_0|_{W_{\lambda}}$, $\im \lambda>0$.  Let
$\lambda=\alpha+i \beta$ with $\beta >0$, then $V_{\lambda}\oplus
V_{\overline{\lambda}}$ is the kernel of the operator $(S_0-\alpha \I)^2+\beta^2 \I$. Define
\[J_{\lambda}= \frac{1}{\beta}(S_0-\alpha \I)|_{W_{\lambda}},\] and
let $(\cdot,\cdot)$ be some scalar product on $W_{\lambda}$.
Note that $J_{\lambda}^2=-\I_{W_{\lambda}}$ and define $(\cdot, \cdot)_{\lambda}$ by
\[
(w,w^{\prime})_{\lambda}:=\frac{1}{2}\left[
  (w,w^{\prime})+(J_{\lambda}w,J_{\lambda}w^{\prime})\right], \quad w,w^{\prime}\in W_{\lambda}.
\]
If $^*$ denotes the involution defined by this scalar product, then
$J_{\lambda}^*J_{\lambda}=\I_{W_{\lambda}}$ i.e., $J_{\lambda}$ is
orthogonal which implies that $J_{\lambda}=-J_{\lambda}^*$.  Also,
$S_0|_{W_{\lambda}}$ is normal since $S_0|_{W_{\lambda}}=\alpha
\I_{W_{\lambda}}+\beta J_{\lambda}$. Moreover, $S_0^*|_{W_{\lambda}}$ maps $W_{\lambda}$ into itself.
\\* \\*
Suppose $\lambda$ is a real eigenvalue of $S_{0}$, then either
$gV_{\lambda} = V_{\lambda}$ or $gV_{\lambda} = V_{\lambda^{-1}}$
depending on whether $\chi(g) = 1$ or $\chi(g) = -1$ respectively.  In
both cases, it is obvious that $gS_{0}^{*}|_{V_{\lambda}} =
gS_{0}|_{V_{\lambda}}
=S^{\chi(g)}_{0}|_{V_{\lambda^{\chi(g)}}}g|_{V_{\lambda}} =
(S^*_{0}|_{V_{\lambda^{\chi(g)}}})^{\chi(g)} g|_{V_{\lambda}}$ and we
are only left to prove this equality for eigenvalues with a strictly
positive imaginary part.  If $\chi(g)=1$, then $gW_{\lambda} = W_{\lambda}$ and therefore $gS_{0}|_{W_{\lambda}}
= S_{0}|_{W_{\lambda}} g|_{W_{\lambda}}$.  But then $g|_{W_{\lambda}}$
commutes with $S_{0}^{*}|_{W_{\lambda}}$ since $S_{0}|_{W_{\lambda}}$
is normal.  In case $\chi(g)=-1$, $gW_{\lambda} = W_{\lambda^{-1}}$,
$J_{\lambda^{-1}} = \frac{1}{\beta}(S_{0}^{-1} - \alpha
I)|_{W_{\lambda^{-1}}}$ and
$$S_{0}|_{W_{\lambda^{-1}}} = \frac{\alpha}{\alpha^2 +
  \beta^2}I|_{W_{\lambda^{-1}}} - \frac{\beta}{\alpha^2 + \beta^2} J_{\lambda^{-1}}.$$
Straightforward calculation reveals that $g J_{\lambda} g^{-1}|_{W_{\lambda^{-1}}} = J_{\lambda^{-1}}$ and therefore
$$ g S_{0}^{*}|_{W_{\lambda}} g^{-1}|_{W_{\lambda^{-1}}} =
 (\alpha I|_{W_{\lambda^{-1}}} - \beta J_{\lambda^{-1}}) = (S_{0}^{*}|_{W_{\lambda^{-1}}})^{-1}. $$
\end{proof}

\begin{lemma}\label{sp}
Let $S_{0}$, $\mathcal{G}$ and $\chi$ be as before, then there exists
a scalar product on $\R^{n}$ (with corresponding involution $*$) such
that $S_{0}S_{0}^{*} = S_{0}^{*}S_{0}$ and $g^{*}g = I$ for all $g \in \mathcal{G}$
\end{lemma}
\begin{proof}
Denote by $\left( \cdot, \cdot \right)$ and $^{T}$ the scalar product
and associated involution given by  Lemma \ref{inprod1}.  The
remainder of the proof then consists in showing that the scalar
product $\left\langle \cdot, \cdot \right\rangle$ defined by
$$ \left\langle x , y \right\rangle = \sum_{g \in \mathcal{G}} \left(
  gx , gy \right) $$ has all the desired properties.  Clearly, $g^{*}g
  = I$ for all $g \in \mathcal{G}$, since the right translation on
  $\mathcal{G}$ defined by $g$ is a bijection.  Therefore, we should
  only prove that $S_{0}^{*}S_{0} = S_{0}S_{0}^{*}$, which is a direct consequence of the following calculation:
\begin{eqnarray*}
\left\langle S_{0}x ,S_{0}y \right\rangle &=& \sum_{g \in \mathcal{G}}
\left( gS_{0}x , gS_{0}y \right) \\
& = & \sum_{g \in \mathcal{G}} \left( gx , (S_{0}^{\chi(g)})^{T} S_{0}^{\chi(g)} gy \right) \\
& = & \sum_{g \in \mathcal{G}} \left( gx , S_{0}^{\chi(g)}
  (S_{0}^{\chi(g)})^{T}gy \right) \\
& = &\sum_{g \in \mathcal{G}} \left( gS_{0}^{T}x , gS_{0}^{T}y \right) \\
& = & \left\langle S_{0}^{T}x ,S_{0}^{T}y \right\rangle \\
\end{eqnarray*}
and the direct observation that $S_{0}^{T} = S_{0}^{*}$.
\end{proof}

\begin{Remark} {\rm{If $\mathcal{G}$ is a compact continuos group with left invariant
  Haar measure $dg$, one may define a $\mathcal{G}$-invariant inner product on
  $\R^n$ by averaging an inner product for $\R^n$ over $\mathcal{G}$:
\[\langle x,y\rangle = \int_\mathcal{G} (hx,hy)dg, \quad x,y\in \R^n
  .\] It can then be shown that every real representation is
  isomorphic to an orthogonal representation. So, one obtains the
  counterparts of Lemma's \ref{inprod1} and \ref{sp} for the continuous case. }}
\end{Remark}

\subsection{Linear Nilpotent Normal Form}
\begin{lemma}
Let $A_0=S_0e^{\mathcal{N}_0}$ be the SU-decomposition of $A_0\in GL_{\mathcal{G}}^{\chi}(n,\R)$ \cite{CMC-Vdb 2}. Then,
$\Ad(A_0^{-1})-\I: U_1\rightarrow V_1$ is an isomorphism as well as $\Ad(A_0^{-1})-\I: U_2\rightarrow V_2$,
where
\begin{eqnarray*}
& & U_1=\im(\Ad(S_0^{-1})-\I)\cap gl_{\mathcal{G}}^{1}(n,\R), \\
& & U_2=\im(\Ad(S_0^{-1})-\I)\cap gl^{\chi}_{\mathcal{G}}(n,\R),\\
& & V_1=\im(\Ad(S_0^{-1})-\I)\cap gl^{\widetilde{\chi}}_{G^{\chi}(A_0)}(n,\R), \\
& & V_2=\im(\Ad(S_0^{-1})-\I)\cap gl^1_{G^{\chi}(A_0)}(n,\R).
\end{eqnarray*}
\end{lemma}

\begin{proof}
We know that $(\Ad(A^{-1}_0)-\I)$ is invertible on
$\im(\Ad(S^{-1}_0)-\I)$ and that it maps  $U_1$ in $V_1$ and $U_2$
into $V_2$ injectively.  We show that $\Ad(A_0^{-1}-\I)$ is an
isomorphism from $U_2$ to $V_2$ by proving surjectivity.  \\*
Let $\psi\in gl^{1}_{\mathcal{G}^{\chi}(A_0)}(n,\R)\cap
\im(\Ad(S^{-1}_0)-\I) $ and let $\widetilde{\psi}\in
\im(\Ad(S^{-1}_0)-\I)$ be such that
$\psi=(\Ad(A_0^{-1})-\I)\widetilde{\psi}$ ($\widetilde{\psi}$ exists
and is unique). Then, one calculates that
\[\psi = P^{1}_{\mathcal{G}^{\chi}(A_0)}\psi=P^{1}_{\mathcal{G}^{\chi}(A_0)}(\Ad(A_0^{-1})-\I)\widetilde{\psi}=
((\Ad(A_0^{-1})-\I))P^{\chi}_{\mathcal{G}}\widetilde{\psi}.
\]
So, if we can show that $P^{\chi}_{\mathcal{G}}\widetilde{\psi}\in
\im(\Ad(S^{-1}_0)-\I)$ then
$P^{\chi}_{\mathcal{G}}\widetilde{\psi}=\widetilde{\psi}$ by
uniqueness and the surjectivity is proved. Now, we know that there
exists some $\widehat{\psi}\in gl(n,\R)$ such that
$\widetilde{\psi}=(\Ad(S_0^{-1})-\I)\widehat{\psi}$. Using the fact
that $S_0^{-1}\in GL_{\mathcal{G}}^{\chi}(n,\R)$ one then obtains that
\[P_{\mathcal{G}}^{\chi}\widetilde{\psi}=
P_{\mathcal{G}}^{\chi}\left((\Ad(S_0^{-1})-\I)\widehat{\psi}\right)=(\Ad(S_0^{-1})-\I)
P^{1}_{\mathcal{G}^{\chi}(S_0^{-1})}\widehat{\Psi},
\]
and the result follows.
\end{proof}

\begin{lemma}\label{LSNF}
Given $A_0\in GL_{\mathcal{G}}^{\chi}(n,\R)$, let
$A_0=S_0e^{\mathcal{N}_0}$ be its SU-decomposition. Then there
exist a neighborhood $\Omega$ of $A_0$ in
$GL_{\mathcal{G}}^{\chi}(n,\R)$ and a map $\phi:\Omega \rightarrow
gl^{1}_{\mathcal{G}}(n,\R) \cap \im(\Ad(S_{0}^{-1}) - \I)$ such
that $\phi(A_0)=0$ and
\[\Ad\left(e^{\phi(A)}\right)A=A_0e^{B(A)},
\]
where $B(A)\in \ker(\Ad(S_0)-\I)\cap gl^{\widetilde{\chi}}_{\mathcal{G}^{\chi}(A_0)}$ and $B(A_0)=0$.
\end{lemma}

\begin{proof}
Define $f:gl_{\mathcal{G}}^1(n,\R)\times GL_{\mathcal{G}}^{\chi}(n,\R)\rightarrow
gl_{\mathcal{G}^{\chi}(A_0)}^{\widetilde{\chi}}(n,\R)$ by the relation
\begin{equation}
\Ad(e^\phi)A=A_0e^{f(\phi,A)}.
\end{equation}
That is, $f(\phi,A):=\log(g(\phi,A))$ with $g(\phi,A):=
A_0^{-1}e^{\phi}Ae^{-\phi}$. The map $f$ is well-defined and smooth for
$(\phi, A)$ near $(0,A_0)$ since $g(0,A_0)=\I$. Also, $\Ad(e^{\phi})A\in GL_{\mathcal{G}}^{\chi}(n,\R )$ and
\begin{eqnarray*}
&& f(0,A_0)=0  \\
&& D_{\phi}f(0,A_0)=(\Ad(A_0^{-1})-\I)|_{gl_{\mathcal{G}}^1(n,\R)}\in
\mathcal{L}(gl_{\mathcal{G}}^1(n,\R), gl^{\widetilde{\chi}}_{\mathcal{G}^\chi(A_0)}(n,\R)).
\end{eqnarray*}
Consider the splitting
\begin{multline}\label{splitting}
gl_{\mathcal{G}^{\chi}(A_0)}^{\widetilde{\chi}}(n,\R)=
\left(gl_{\mathcal{G}^{\chi}(A_0)}^{\widetilde{\chi}}(n,\R)\cap\im(\Ad(S_0^{-1})-\I)\right)\\
\oplus\left(gl_{\mathcal{G}^{\chi}(A_0)}^{\widetilde{\chi}}(n,\R)\cap\ker(\Ad(S_0^{-1})-\I)\right)
\end{multline}
and let
$\pi:gl_{\mathcal{G}^{\chi}(A_0)}^{\widetilde{\chi}}(n,\R)\rightarrow
\left(
  gl_{\mathcal{G}^{\chi}(A_0)}^{\widetilde{\chi}}(n,\R)\cap\im(\Ad(S_0^{-1})-\I)\right)$
the corresponding projection on the first factor. Define
\[h:=\pi\cdot f|_{gl^{1}_{\mathcal{G}}(n,\R)\times GL_{\mathcal{G}}^\chi(n,\R)}.\]
Then, $h(0,A_0)=0$ and the operator $D_{\phi}h(0,A_0)$ is an
isomorphism between
$gl_{\mathcal{G}^{\chi}(A_0)}^{\widetilde{\chi}}(n,\R)\cap\im(\Ad(S_0^{-1})-\I)$
and
$gl_{\mathcal{G}}^{1}(n,\R)\cap\im(\Ad(S_0^{-1})-\I)$.
The Implicit Function Theorem then implies that there exists
$\widetilde{\phi}: GL_{\mathcal{G}}^{\chi}(n,\R)\rightarrow
gl^1_{\mathcal{G}}(n,\R)$
with $\widetilde{\phi}(A_0)=0$ such that
$h(\widetilde{\phi}(A),A)=0$. Hence, setting
$B(A):=f(\widetilde{\phi}(A),A)\in \ker(\Ad(S_0)-\I)\cap gl_{\mathcal{G}^{\chi}(A_0)}^{\widetilde{\chi}}(n,\R)$
proves the lemma.
\end{proof}

\noindent
{\bf{Proof of the splitting (\ref{splitting}).}}
We prove that the splitting (\ref{splitting}) holds.
The semisimplicity of $S_0\in GL_{\mathcal{G}}^{\chi}(n,\R)$
implies that
$gl(n,\R)=\im(\Ad(S_0^{-1})-\I)\oplus\ker(\Ad(S_0^{-1})-\I)$.
Therefore, each $\Psi\in gl(n,\R)$ has unique decomposition
\begin{equation}\label{split_Psi}
\Psi=\Psi_{im}+\Psi_{ker}
\end{equation}
 with
$\Psi_{im}\in \im(\Ad(S_0^{-1})-\I)$ and $\Psi_{ker} \in \ker(\Ad(S_0^{-1})-\I)$.
Recall that $\ker(\Ad(S_0^{-1})-\I)=\ker(\Ad(S_0)-\I)$.
For $\Psi\in gl^{\widetilde{\chi}}_{\mathcal{G}^{\chi}(A_0)}(n,\R)$ this allows us to write:
$$\Psi=P^{\widetilde{\chi}}_{\mathcal{G}^{\chi}(A_0)}\Psi=P^{\widetilde{\chi}}_{\mathcal{G}^{\chi}(A_0)}
\Psi_{im}+P^{\widetilde{\chi}}_{\mathcal{G}^{\chi}(A_0)}\Psi_{ker}.$$
Using the fact that $S_0^{-1}\in GL_{\mathcal{G}}^{\chi}(n,\R)$ and the definitions, one checks that
\begin{eqnarray*}
&& \left(\Ad(S_0)-\I\right)P^{\widetilde{\chi}}_{\mathcal{G}^{\chi}(A_0)}=P^{1}_{\mathcal{G}^{\chi}
(A_0S_0^{-1})} \left(\Ad(S_0)-\I\right) = P^{\chi}_{\mathcal{G}^{\chi}(A_{0})} (\Ad(S_0^{-1}) - \I), \\
&& P^{\widetilde{\chi}}_{\mathcal{G}^{\chi}(A_0)}\left(\Ad(S_0^{-1})-\I\right)=
\left(\Ad(S_0^{-1})-\I\right)P^{1}_{\mathcal{G}^{\chi}(A_0S_0^{-1})}.
\end{eqnarray*}
Therefore, $P^{\widetilde{\chi}}_{\mathcal{G}^{\chi}(A_0)}\Psi_{im}\in \im(\Ad(S_0^{-1})-\I)$ and
$P^{\widetilde{\chi}}_{\mathcal{G}^{\chi}(A_0)}\Psi_{ker}\in \ker(\Ad(S_0)-\I)$.
The uniqueness of the splitting (\ref{split_Psi}) then implies that
$\Psi_{im}=P^{\widetilde{\chi}}_{\mathcal{G}^{\chi}(A_0)}\Psi_{im}$
and
$\Psi_{ker}=P^{\widetilde{\chi}}_{\mathcal{G}^{\chi}(A_0)}\Psi_{ker}$. Hence,
$\Psi_{im}\in \im(\Ad(S_0^{-1})-\I)\cap
gl^{\widetilde{\chi}}_{\mathcal{G}^{\chi}(A_0)}(n,\R)$ and
$\Psi_{ker}\in \ker(\Ad(S_0^{-1})-\I)\cap
gl^{\widetilde{\chi}}_{\mathcal{G}^{\chi}(A_0)}(n,\R)$,
which proves the result.

\begin{proposition}[Linear Nilpotent Normal Form]\label{LNNF}
Let $A_0\in GL_{\mathcal{G}}^{\chi}(n,\R)$ have SU-decomposition
$A_0=S_0e^{\mathcal{N}_0}$ and let $\left\langle \cdot , \cdot
\right\rangle$ be a scalar product as in Lemma $4$.  Then, there
exist a neighborhood $\Omega$ of $A_0$ in
$GL_{\mathcal{G}}^{\chi}(n,\R)$ and a map $\phi:\Omega \rightarrow
gl^{1}_{\mathcal{G}}(n,\R)$ such that $\phi(A_0)=0$ and
\[
\Ad\left(e^{\phi(A)}\right)A=S_0e^{\mathcal{N}_0+C(A)},
\]
where $C(A)\in \ker(\Ad(S_0)-\I)\cap \ker(\ad(\mathcal{N}_0^T))\cap gl^{\chi}_{\mathcal{G}}(n,\R)$ and $C(A_0)=0$.
\end{proposition}

\begin{proof}
From Lemma \ref{LSNF} we can assume that there exists a neighborhood
$\Omega^{*} \subset GL^{\chi}_{\mathcal{G}} (n, \R)$ of $A_{0}$ and a
map $\widetilde{\phi} : \Omega^{*} \rightarrow gl^{1}_{\mathcal{G}}(n , \R)$ such that
\[\Ad(e^{\widetilde{\phi}(A)})A=S_0e^{\mathcal{N}_0+B(A)},\]
where $B(A)\in \ker(\Ad(S_0)-\I)\cap gl^{\chi}_{\mathcal{G}}(n,\R)$ and $B(A_0)=0$.
Let $\phi \in \ker(\Ad(S_0)-\I) \cap gl^{1}_{\mathcal{G}}(n,\R)$, then,
\[\Ad(e^{\phi}) \Ad(e^{\widetilde{\phi}(A)}) A=S_0e^{\mathcal{N}_0+f(\phi,A)},\]
where
$f:\left(\ker(\Ad(S_0)-\I)\cap gl^{1}_{\mathcal{G}}(n,\R)\right)
\times GL^{\chi}_{\mathcal{G}}(n,\R)\rightarrow gl^{\chi}_{\mathcal{G}}(n,\R)$ is given by
\[f(\phi,A):=\Ad(e^{\phi})\left(\mathcal{N}_0+B(A)\right)-\mathcal{N}_0.\]
We have that $f(0,A_0)=0$ and that
$D_{\phi}f(0,A_0)=-\ad(\mathcal{N}_0)|_{\ker(\Ad(S_0)-\I)\cap gl^{1}_{\mathcal{G}}(n,\R)}$
belongs to $\mathcal{L}\left(\ker(\Ad(S_0)-\I)\cap
  gl^{1}_{\mathcal{G}}(n,\R), \ker(\Ad(S_0)-\I)
\cap gl^{\chi}_{\mathcal{G}}(n,\R)\right).$
Now, by Lemma $4$, the operators $\ad(\mathcal{N}_0)$ and
$\ad(\mathcal{N}_0^T)$
leave the complementary subspaces $\ker(\Ad(S_0)-\I)$ and $\im(\Ad(S_0)-\I)$ of $gl(n,\R)$ invariant and
\[\left(\ad(\mathcal{N}_0)|_{\ker(\Ad(S_0)-\I)}\right)^T=\ad(\mathcal{N}_0^T)|_{\ker(\Ad(S_0)-\I)}.\]
Also, $\ad(\mathcal{N}_0)$ maps $gl^1_{\mathcal{G}}(n,\R)$ on $gl^\chi_{\mathcal{G}}(n,\R)$ and viceversa. Hence,
\begin{multline}\label{split_1}
\ker\left(\Ad(S_0)-\I\right)\cap gl^\chi_{\mathcal{G}}(n,\R)= \left[
\ker\left(\Ad(S_0)-\I\right)\cap \im(\ad(\mathcal{N}_0))\cap gl^\chi_{\mathcal{G}}(n,\R)\right] \\
\oplus \left[\ker\left(\Ad(S_0)-\I\right)\cap\ker(\ad(\mathcal{N}_0^T))\cap gl^\chi_{\mathcal{G}}(n,\R)\right]
\end{multline}
and
\begin{multline}\label{split_2}
\ker\left(\Ad(S_0)-\I\right)\cap gl^1_{\mathcal{G}}(n,\R)= \left[
\ker\left(\Ad(S_0)-\I\right)\cap \ker(\ad(\mathcal{N}_0))\cap gl^1_{\mathcal{G}}(n,\R)\right] \\
\oplus \left[\ker\left(\Ad(S_0)-\I\right)\cap\im(\ad(\mathcal{N}_0^T))\cap gl^1_{\mathcal{G}}(n,\R)\right] .
\end{multline}
Note that (by the choice of the scalar product) the operator
$\ad(\mathcal{N}_0)$ is an
isomorphism from
$\left[\ker\left(\Ad(S_0)-\I\right)\cap\im(\ad(\mathcal{N}_0^T))\cap
  gl^\chi_{\mathcal{G}}(n,\R)\right]$ onto
$\left[\ker\left(\Ad(S_0)-\I\right)\cap \im(\ad(\mathcal{N}_0))\cap
  gl^\chi_{\mathcal{G}}(n,\R)\right]$. Indeed,
$\ad(\mathcal{N}_{0}) P^1_{\mathcal{G}} = P^{\chi}_{\mathcal{G}} \ad(\mathcal{N}_{0})$.
Let $$\pi: \ker\left(\Ad(S_0)-\I\right)\cap gl^\chi_{\mathcal{G}}(n,\R)\rightarrow \left[
\ker\left(\Ad(S_0)-\I\right)\cap \im(\ad(\mathcal{N}_0))\cap
gl^\chi_{\mathcal{G}}(n,\R)\right]$$
be the projection on the first factor associated to the splitting (\ref{split_1}), and define the map
\begin{multline*}
g: \left(
  \ker\left(\Ad(S_0)-\I\right)\cap\im(\ad(\mathcal{N}_0^T))\cap
  gl^1_{\mathcal{G}}(n,\R) \right)
\times GL^{\chi}_{\mathcal{G}} (n, \R)\\
\rightarrow \ker\left(\Ad(S_0)-\I\right)\cap \im(\ad(\mathcal{N}_0))\cap gl^\chi_{\mathcal{G}}(n,\R),
\end{multline*}
by
\[g:=\pi \cdot f|_{\left(
    \ker\left(\Ad(S_0)-\I\right)\cap\im(\ad(\mathcal{N}_0^T))
\cap gl^1_{\mathcal{G}}(n,\R) \right) \times GL^{\chi}_{\mathcal{G}} (n, \R)}.\]
Then, $g(0,A_0)=0$ and $D_{\phi}g(0,A_0)$ is an isomorphism. The
Implicit Function
Theorem then implies the existence of a map $\phi^*:\Omega \subset
\Omega^{*} \rightarrow \ker\left(\Ad(S_0)-\I\right)
\cap\im(\ad(\mathcal{N}_0^T))\cap gl^1_{\mathcal{G}}(n,\R)$ with
$\phi^*(A_0)=0$ such that
\[g(\phi^*(A),A)=0.\] Then,
$f(\phi^*(A),A)\in
\ker\left(\Ad(S_0)-\I\right)\cap\ker(\ad(\mathcal{N}_0^T))\cap
gl^\chi_{\mathcal{G}}(n,\R)$.
Hence, the result follows.
\end{proof}

\section{Normal Form}\label{sec-normalform}
Goal of this section is to prove by induction that a map
$\psi_{\lambda}$ as considered in this paper admits a
$\chi$-equivariant normal form with constraints involving the
nilpotent $\mathcal{N}_0$. More in details, we first prove that a
normalisation as in (\ref{ReductionNFformula1}) is possible. We fix then a scalar
product in $\R^n$ as given in Lemma \ref{sp} and show that
$\psi_{\lambda}^{NF}$ can be brought into the form
$\psi_{\lambda}^{NF}=S_0 \exp(\mathcal{N}_0+X_{\lambda})$ with
$\mathcal{N}_0+X_{\lambda}\in \ker(\Ad(S_0)-\I)\cap
\im(P_{\mathcal{G}^{\chi}})$ and
$DX_{\lambda}(0)\mathcal{N}_0^T(x)=\mathcal{N}_0^TX_{\lambda}(x)$.

\medskip
Let $\mathcal{H}_{k}=\mathcal{H}_k(\R^n)$ be the space of polynomial maps homogeneous of degree
$k$, note that $\mathcal{H}_1=gl(n,\R)$.
Then the Taylor series of $\Phi\in {\rm{Diff}}_0(\R^n)$ is an element of the space of
formal power series $\prod_{k=1}^{\infty}\mathcal{H}_k$. Also, considering
 $A_0=D\Phi_0(0)\in \mathcal{L}(\R^n)$, it directly follows from the definitions that
$\Ad(A_0)$ induces a linear map $\mathcal{H}_k \rightarrow \mathcal{H}_k$ to be denoted
by $\Ad_{k}A_0.$  Let $\mathcal{X}_{0}$ be the space of all smooth
vectorfields on $\R^{n}$ with a fixed point at $0$, and let
$\mathcal{X}_{0}^{k} := \left\{ X \in \mathcal{X}_{0} | D_{x}^{j} X(0)
  = 0, 1 \leq j \leq k \right\}.$
The following notation is self explanatory:
\begin{equation*}
X=X_{1}+\cdots +X_{k}\;\, {\rm{mod}}\mathcal{X}_{0}^{k},\qquad \text{with }X_{j}\in \mathcal{H}_{j}
\label{taylor X v.f.}.
\end{equation*}
Define the operator $C_k: \mathcal{H}_1 \rightarrow
\mathcal{L}\left(\mathcal{H}_k\right)$, for all $X_1\in \mathcal{H}_1$, $X_k\in \mathcal{H}_k$, by
\begin{equation}\label{definition CkX_1}
 C_{k}\left( X_{1}\right) X_{k}:=
\int\limits_{0}^{1}e^{ -s X_{1}}X_{k}e^{sX_{1}} ds.
\end{equation}
Observe also that if $X_1=0$, then
\begin{equation*} C_{k}\left( 0\right) =\I_{\mathcal{H}_{k}}. \end{equation*}
Let $\mathcal{P}_{k}$ be the space of all polynomial maps $p:\R^n \rightarrow \R^n$, $p(0)=0$, of
degree less or equal to $k$. According to our previous notation,
$\mathcal{P}_k=\mathcal{X}_0/\mathcal{X}_0^k$, $k\geq 1$.  Set
\begin{equation*}
X^{\left[ k\right] }:=X_{1}+\cdots +X_{k}\in \mathcal{P}_{k}, \quad \text{with each }X_{j}\in \mathcal{H}_{j}.
\end{equation*} \\*
We need the following generalisation of the Campbell-Hausdorff formula \cite{CMC1}.
\begin{lemma}[\cite{CMCPhD}]
\label{CampbellHausdorf}
Given $X\in {\mathcal{X}_0}$, let $X_1:=X \, {\rm{mod}}\mathcal{X}_0^1$
such that $C_k(X_1)$ is invertible. Then, for any $Y_{k}\in \mathcal{H}_{k}$,
\begin{description}
\item[(i)] $e^{X}e^{Y_{k}}=e^{X+C_{k}\left( X_{1}\right)
^{-1}Y_{k}}\; {\rm{mod}}\mathcal{X}_0^{k} ,$ \\
\item[(ii)] $e^{Y_{k}}e^{X}=e^{X+C_{k}\left( -X_{1}\right)
^{-1}Y_{k}}\; {\rm{mod}}\mathcal{X}_0^{k}.$
\end{description}
\end{lemma}

\begin{proposition} \label{SNF}
Let $\mathcal{G}$ be a finite compact group with one-dimensional
character $\chi:\mathcal{G}\rightarrow \C$. Let $\psi:\R^n\times
\R^m \rightarrow \R^n$ be a smooth family of local diffeomorphisms
satisfying (\ref{symmetry}) and (\ref{basic fixed point}). Let
$A_0=S_0e^{\mathcal{N}_0}\in GL_{\mathcal{G}}^{\chi}(n,\R)$ be the
SU-decomposition of $A_0=D_x\psi_0(0)$. Then, for each $k\geq 1$
there exists a neighborhood $\omega_k$ of the origin in $\R^m$ and
a parameter dependent near-identity transformation
$\Phi_{k,\lambda}:\R^n\rightarrow \R^n$ with
$\Phi_{k,\lambda}(0)=0$, $D_x\Phi_{0}(0)=\I$ and
$g\Phi_{k,\lambda}g^{-1}=\Phi_{k,\lambda}$ such that the following
holds:
\begin{equation}
\Ad(\Phi_{k,\lambda})\psi_{\lambda}=A_0e^{X_{\lambda}^{[k]}} \mod \mathcal{X}_0^k, \quad \forall \lambda \in \omega_k,
\end{equation}
with $X_{\lambda}^{[k]}\in \ker(\Ad(S_0)-\I)\cap
\im(P^{\widetilde{\chi}}_{\mathcal{G}^{\chi}(A_0)})\subset
\mathcal{P}_k$,
and $X^{[k]}_{0}=0\mod \mathcal{X}_0^{1}$.
\end{proposition}

\begin{proof}
The proof is by induction on $k$. For $k=1$ the result follows from
Lemma \ref{LSNF}. The induction step is proved as follows.
Assume the result true for $(k-1)$ ($k>1$). Denoting
$\Ad(\Psi_{k-1,\lambda})(\psi_{\lambda})$ again by $\psi_{\lambda}$, this means that
\[\psi_{\lambda}=A_0e^{X^{[k-1]}_{\lambda}}\mod\mathcal{X}_0^{k-1},\]
with
$X_{\lambda}^{[k-1]}\in
\ker(\Ad(S_0)-\I)\cap\im(P^{\widetilde{\chi}}_{\mathcal{G}^{\chi}(A_0)})
\subset \mathcal{P}_{k-1}$ and $X_{0}^{[k-1]}=0\mod
\mathcal{X}_{0}^{1}$.
Then,
\[\psi_{\lambda}=A_0e^{X^{[k-1]}_{\lambda}+Z_{k,\lambda}}\mod\mathcal{X}_0^{k},\]
for some
$Z_{k,\lambda}\in \mathcal{H}_k\cap
\im(P^{\widetilde{\chi}}_{\mathcal{G}^{\chi}(A_0)})$. The goal is to
show that we can transform $\psi_{\lambda}$ further so that
$Z_{k,\lambda}\in \ker(\Ad(S_0)-\I)\cap
\im(P^{\widetilde{\chi}}_{\mathcal{G}^{\chi}(A_0)})$.
To do so, let $\Phi_{k}$ be of the form $\Phi_{k}=e^{\phi_{k}}$
with $\phi_{k}\in \mathcal{H}_k\cap \im(P^1_{\mathcal{G}})$. Then,
$\Ad(e^{\phi_k})\psi_{\lambda}\in \im(P^{\chi}_{\mathcal{G}})$ and
\[\Ad(e^{\phi_k})\psi_{\lambda}=A_0e^{X^{[k-1]}_{\lambda}+f(\phi_k,\lambda)}\mod \mathcal{X}_0^{k},\]
where
$f_{k}:\mathcal{H}_{k}\cap\im(P^1_{\mathcal{G}})\rightarrow
\mathcal{H}_{k}\cap \im(P^{\widetilde{\chi}}_{\mathcal{G}^{\chi}(A_0)})$ is given by
\[f_k(\phi_k,\lambda)=Z_{k,\lambda}+C_k(-X_{1,\lambda})\Ad(A_0^{-1})\phi_{k}-C_{k}(X_{1,\lambda})^{-1}\phi_k.\]
Note that $f_{k}(\phi_k,\lambda)$ is smooth and well defined near
$(0,0)$ and $f_k(0,\lambda)=Z_{k,\lambda}.$ Using the fact that
$\Ad_{k}(A_0^{-1})-\I$ is an isomorphism between $\mathcal{H}_k\cap\im(P^1_{\mathcal{G}})\cap \im(\Ad_k(S_0^{-1})-\I)$
and $\mathcal{H}_k\cap \im(P^{\widetilde{\chi}}_{\mathcal{G}^\chi(A_0)})\cap \im(\Ad_k(S_0^{-1})-\I)$ and that
\begin{multline*}
\mathcal{H}_k\cap\im(P^{\widetilde{\chi}}_{\mathcal{G}^\chi(A_0)})=\mathcal{H}_k\cap
\im(P^{\widetilde{\chi}}_{\mathcal{G}^\chi(A_0)})\cap \im(\Ad_k(S_0^{-1})-\I)\\
\oplus
\mathcal{H}_k\cap \im(P^{\widetilde{\chi}}_{\mathcal{G}^\chi(A_0)})\cap \ker(\Ad_k(S_0^{-1})-\I),
\end{multline*}
the result follows by the Implicit Function Theorem as in Lemma \ref{LSNF}.
\end{proof}

\begin{theorem}[Nilpotent Normal Form]\label{NNF}
Let $\mathcal{G}$ be a finite compact group with one-dimensional
character $\chi:\mathcal{G}\rightarrow \C$. Let $\psi:\R^n\times
\R^m \rightarrow \R^n$ be a smooth family of local diffeomorphisms
satisfying (\ref{symmetry}) and (\ref{basic fixed point}). Let
$A_0=S_0e^{\mathcal{N}_0}\in GL_{\mathcal{G}}^{\chi}(n,\R)$ be the
SU-decomposition of $A_0$, and let $\langle\cdot, \cdot \rangle$
be a scalar product as in Lemma \ref{sp}. Then, for each $k\geq 1$
there exists a neighborhood $\omega_k$ of the origin in $\R^m$ and
a parameter dependent near-identity transformation
$\Phi_{k,\lambda}:\R^n\rightarrow \R^n$ with
$\Phi_{k,\lambda}(0)=0$, $D_x\Phi_{0}(0)=\I$ and
$g\Phi_{k,\lambda}g^{-1}=\Phi_{k,\lambda}$ such that the following
holds:
\begin{equation}
\Ad(\Phi_{k,\lambda})\psi_{\lambda}=S_0e^{\mathcal{N}_0+X_{\lambda}}\mod\mathcal{X}_0^k,
\quad \forall \lambda \in \omega_k,
\end{equation}
with
\begin{eqnarray}\label{propvf}
& & X_{\lambda}(0)=0 \quad DX_{0}(0)=0, \notag\\
& & S_0\cdot X_{\lambda}=X_{\lambda}\cdot S_0,\notag\\
& & DX_{\lambda}(x)\mathcal{N}_0^T(x)=\mathcal{N}_0^TX_{\lambda}(x), \\
& & (\mathcal{N}_0+X_{\lambda})\cdot g=\chi(g)g\cdot (\mathcal{N}_0+X_{\lambda})\notag.
\end{eqnarray}
\end{theorem}
\begin{proof}
By Proposition \ref{SNF} we can assume that after transformation $\psi_{\lambda}$ is of the form
\[\psi_{\lambda}=S_0e^{\mathcal{N}_0+X^{[k]}_{\lambda}}\mod \mathcal{X}_0^k,\] where
$X^{[k]}_{\lambda}\in \ker(\Ad(S_0)-\I)\cap\im(P^{\chi}_{\mathcal{G}})$ and
$X^{[k]}_{\lambda}(0)=0$, $DX^{[k]}_{0}(0)=0$.
The proof proceeds then by induction on $k$. For $k=1$, the result
follows from Proposition \ref{LNNF}. For $k>1$, the induction argument
is similar to that of Proposition \ref{SNF} with the following change:
\[\psi_{\lambda}=S_0e^{\mathcal{N}_0+X^{[k-1]}_{\lambda}}\mod\mathcal{X}_0^{k-1},\]
with $X_{\lambda}^{[k-1]}\in \ker(\Ad(S_0)-\I)\cap\ker(\ad(\mathcal{N}_0^T))\cap\im(P^{\chi}_{\mathcal{G}})\subset
\mathcal{P}_{k-1}$ and $X_{0}^{[k-1]}=0\mod \mathcal{X}_{0}^{1}$.
Then,
\[\psi_{\lambda}=A_0e^{\mathcal{N}_0+X^{[k-1]}_{\lambda}+Z_{k,\lambda}}\mod\mathcal{X}_0^{k},\]
for some $Z_{k,\lambda}\in\ker(\Ad_k(S_0)-\I)\cap\im(P^{\chi}_{\mathcal{G}})$. The goal is to
bring the term $Z_{k,\lambda}$ in $\ker(\Ad(S_0)-\I)\cap\ker(\ad(\mathcal{N}_0^T))\cap\im(P^{\chi}_{\mathcal{G}})$
which is achieved by a similar splitting argument as in Proposition \ref{LNNF}.
\end{proof}

\section{Example}
Consider the involution $R:=\left(\begin{array}{cc} 0& 1 \\ 1& 0
  \end{array}\right)\in \mathcal{L}(\R^2)$ and the finite
compact group $\mathcal{G}=\{\I,R\}$ with character
$\chi:\mathcal{G}\rightarrow \R$ given by $\chi(I)=1$
and $\chi(R)=-1$. Let $\Phi_1:\Omega \subset \R^2\rightarrow\R^2$ be the map given by
\[(x,y)\mapsto (\frac{x^3}{y^2},\frac{x^2}{y}),\] here
$\Omega:=\{(x,y)\in \R^2|\, x>0,\ y>0\}$ is the open first
quadrant. Now, $\Phi_1$ is a local $R$-reversible diffeomorphism,
that is, $R\cdot \Phi_{1}(x,y)=\Phi_{1}^{-1}(R(x,y)).$ Note that
the map $\Phi_1$ has a line of fixed points, $y=x$, that is
invariant under $R$. Indeed, $R(x,x)=(x,x)$. We consider the map
$\Phi_1$ around the fixed point $P:=(1,1)$ and calculate its
normal form (up to the third order) using the theory developed in
this paper. The linearisation of $\Phi_1$ at the fixed point $P$
is given by
\[A_0:=\left(\begin{array}{cc}3& -2 \\ 2 & -1 \end{array}\right)=S_0+N_0,\]
where $S_0=\I$ and $N_0:=\left(\begin{array}{cc}2& -2 \\ 2 & -2 \end{array}\right)$.
Theorem \ref{NNF} tells us that there exists an $R$-equivariant (and $S_0$-equivariant) transformation
$g\in \mathcal{H}_2$ such that
\begin{equation}\label{nfex1}\Ad(e^g)\Phi=S_0e^{\mathcal{N}_0+X_{2}}\mod \mathcal{X}^3_0,
\end{equation}
where
$\Phi(x,y)=\Phi_1(1+x,1+y)-(1,1)$ and $X_2\in \mathcal{P}_2$ satisfies the properties (\ref{propvf}).
Recall that $\mathcal{N}_0=\log(\I+S_0^{-1}N_0)$.
Direct calculations show that admissible $X_2$'s must have the form
$$X_2(x,y)=\left( d(x+y)^2,-d(x+y)^2\right),$$
for some $d\in \R$, while $g(x,y)\in \mathcal{H}_2$ must have the form
$$g(x,y)=\left(ax^2+bxy+cy^2, cx^2+bxy+ay^2 \right),$$ for some
$a,b,c\in \R$. So, the relation (\ref{nfex1})
implies that $b = -2c,\  a = -1/2+c,\ d = 0, \ c = c$. That is, there
exists a one-parameter family of maps
$g(x,y)=((c-1/2)x^2-2cxy+cy^2, cx^2-2cxy+(c-1/2)y^2)$ that brings
$\Phi$ into the (third order) normal form
with $X_2=0$. Choosing, for example, $c=0$ one has that
\[\Ad(\exp(-\frac{1}{2}x^2,\frac{1}{2}y^2))\Phi(x,y)=A_0
\left(\begin{array}{c}x \\y\end{array}\right) \mod \mathcal{X}_0^3.\]

\smallskip
\noindent
\begin{Remark}
{\rm{
Note that the map $\Phi_1(x,y)=(\frac{x^3}{y^2},\frac{x^2}{y})$ can be
written as the composition
$\Phi_1(x,y)=R\cdot\Psi(x,y)$ where $\Psi: \Omega\rightarrow \R^2$ is
given by
$\Psi(x,y)=(\Psi_1(x,y),\Psi_2(x,y))=(\frac{x^2}{y},\frac{x^3}{y^2})$.
We show that the map $\Psi$ is not the time-one
map of vector field. Observe that this doesn't say anything about the
map $\Phi_1$.
The map $\Psi$ has the same line of fixed points as $\Phi_1$; i.e.,
$\Psi(x,x)=(x,x)$ and observe
that $\Psi_2(x,y)/\Psi_1(x,y)=x/y$. In polar coordinates we have that
\[\Psi(r,\theta)=(r{\rm{cotg}}^2(\theta),\pi/2-\theta),\]
which shows that points in the first quadrant under the first
bisectrice $y=x$ are mapped into points above this
line as shown in Fig.\ \ref{figura}.

\def\Red#1{{\color{red}#1}}
\def\Green#1{{\color{darkgreen}#1}}
\def\Blue#1{{\color{blue}#1}}

\definecolor{red}{rgb}{1,0,0}
\definecolor{blue}{rgb}{0,0,1}
\definecolor{green}{rgb}{0,1,0}
\definecolor{black}{rgb}{0,0,0}
\definecolor{yellow}{rgb}{1,1,0}
\definecolor{mdwblue}{rgb}{0.2,0.2,0.6}
\definecolor{gray}{rgb}{0.7,0.7,0.7}
\definecolor{darkgreen}{rgb}{0.2,0.7,0.2}

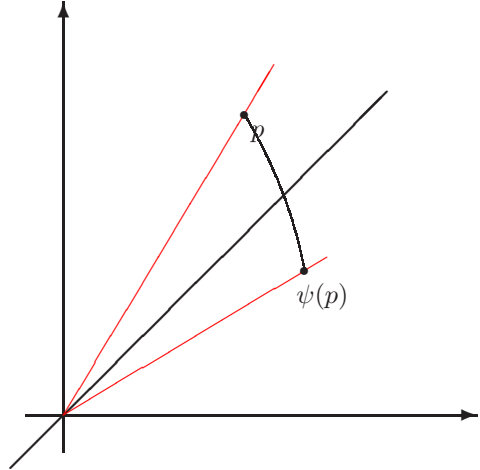
\begin{figure}[t]
\begin{center}
  \setlength{\unitlength}{1cm}
\begin{picture}(8,5)
\thicklines
\put(0.5,1){\vector(1,0){6}}
\put(1,0.5){\vector(0,1){6}}
\put(0.3,0.3){\line(1,1){5}}
\thinlines
\put(1,1){\Red{\line(3,5){2.8}}}
\put(1,1){\Red{\line(5,3){3.5}}}

\put(3.4,5){\circle*{0.1}}
\put(3.5,4.7){$p$}
\put(4.2,2.92){\circle*{0.1}}
\put(4.1,2.5){$\psi(p)$}
\qbezier(3.4,5)(4,4)(4.2,2.92)
\put(4.2,2.92){\line(1,2){0.1}}
\put(4.2,2.92){\line(-2,1){0.1}}

\end{picture}
\caption{The map $\Psi$ maps points under the bisectrice $y=x$ to
  points above it.}
\label{figura}
\end{center}
\end{figure}

Now, suppose that there exists a vector field $X$ which is defined
in a neighborhood of $(1,1)$ and is such that
$\Psi(x,y)=e^{X}(x,y)$. Then, one sees that an integral curve of
$X$ cannot cross the fixed line $y=x$ since otherwise one can find
points which remain in the same half plane defined by the first
bisectrice. Hence, the only remaining possibility is that
$X(x,x)=0$ which is also impossible by the same argument as above.
}}
\end{Remark}

\paragraph{Acknowledgment} The authors wish to thank Prof.\ A.\
Vanderbauwhede, Prof.\ M. Field,  Dr.\ J. Lamb and B.\
Malengier for useful discussions.

\end{document}